\documentstyle[12pt]{amsart}
\def\myperemail{{phung@@ioit.ncst.ac.vn}}

\def\myperaddress{Hanoi Institute of Mathematics, P.O.Box 631, 10000 Boho, Hanoi}

\def\myauthor{Ph\`ung H{{\accent"5E o}\kern-.28em\raise.2ex\hbox{\char'22}\kern-.20em} H{a\kern-.370em\raise.16ex\hbox{\char'47}\kern.1em}i}
\def\myAUTHOR{ PH\`UNG H{{\accent"5E O}\kern-.38em\raise.8ex\hbox{\char'22}\kern-.12em}  H{A\kern-.46em\raise.80ex\hbox{\char'47}\kern.18em}I}

\def\amshead{
\title[Integral on Hopf algebras ]{Integrals on Hopf algebras and Application to Representation Theory of Quantum Groups of Type $A_{0|0}$}
\author{ PH\`UNG H{{\accent"5E O}\kern-.38em\raise.8ex\hbox{\char'22}\kern-.12em}  H{A\kern-.46em\raise.80ex\hbox{\char'47}\kern.18em}i}
\address{\myperaddress}
\keywords{}
\email{\myperemail}
\maketitle }

\def\lora{\longrightarrow}
\def\ot{\otimes}

\def\loma{\longmapsto}

\def\alph{\alpha}
\def\si{\sigma}
\newcommand{\bbas}{\begin{eqnarray*}}
\newcommand{\eeas}{\end{eqnarray*}}
\newcommand{\bbar}{\begin{array}}
\newcommand{\eear}{\end{array}}
\newcommand{\bbs}{\begin{displaymath}}
\newcommand{\ees}{\end{displaymath}}
\newcommand{\bb}{\begin{equation}}
\newcommand{\eqbb}{\begin{equation}}
\def\ee{\end{equation}}
\def\eqee{\end{equation}}
\def\eea{\end{eqnarray}}

\def\bba{\begin{eqnarray}}
\newtheorem{thm}{Theorem}[section]
\newtheorem{lem}[thm]{Lemma}

\newtheorem{rem}[thm]{Remark}

\newtheorem{cor}[thm]{Corollary}

\newtheorem{pro}[thm]{Proposition}

\def\Ker{\mbox{\rm Ker}}
\def\Im{\mbox{\rm Im}}

\def\Hom{\mbox{\rm Hom}}

\def\End{\mbox{\rm End}}

\def\H{{\cal H}}

\def\A{{\cal A}}

\def\L{{\cal L}}
\def\N{{\cal N}}

\def\P{{\cal P}}
\def\D{{\cal D}}

\def\eee{\rule{.75ex}{1.5ex}\\[1ex]}

\def\proof{{\it Proof.\ }}
\newcommand{\va}{\varepsilon}

\newcommand{\rk}{\mbox{\rm rank}}
\def\rank{\rk}

\renewcommand{\dim}{\mbox{\rm dim}}

\def\rref#1{(\ref{#1})}
\newcommand{\bK}{{\mathbf k}}

\def\db{{\mathchoice{\mbox{\rm db}}
                    {\mbox{\rm db}}
                    {\mbox{\scriptsize\rm db}}
                    {\mbox{\tiny\rm db}} }}

\def\ev{{\mathchoice{\mbox{\rm ev}}
                    {\mbox{\rm ev}}
                    {\mbox{\scriptsize\rm ev}}
                    {\mbox{\tiny\rm ev}} }}

\def\id{{\mathchoice{\mbox{\rm id}}
                    {\mbox{\rm id}}
                    {\mbox{\scriptsize\rm id}}
                    {\mbox{\tiny\rm id}} }}

\def\op{{\mathchoice{\mbox{\rm op}}
                    {\mbox{\rm op}}
                    {\mbox{\scriptsize\rm op}}
                    {\mbox{\tiny\rm op}} }}

\newcommand{\Mat}{\mbox{\rm M}}

\def\part{\vdash}

\def\lam{{\lambda}}

\def\I{{\cal I}}

\let\prob=\int
\def\int{\displaystyle\prob}
\newtheorem{proc}[thm]{}
\def\intl#1{\int_l\left({#1}\right)}
\def\intr#1{\int_r\left({#1}\right)}
\def\Cf{{\cal C}{\it f}}
\def\J{{\cal J}}
\def\leftact{\rightharpoonup}
\def\rightact{\leftharpoonup}

\def\sora{\rightarrow}

\def\op{\oplus}

\def\I#1#2{{{\rm I}_{#1,#2}}}

\begin{document}
\bibliographystyle{plain}
\amshead
\section*{Introduction}

In this work we study some properties of (non-cosemisimple) Hopf algebras, possessing integrals, which are also called co-Frobenius Hopf algebras. We apply the result obtained to the classification of representations of quantum groups of type $A_{0|0}$.

The notion of integral on Hopf algebras is motivated by the Haar integral on compact groups. In fact, the axiom of the Haar integral on a compact group can be given in a pure algebraic way as a linear functional on the algebra of (regular) function on the group (which is a Hopf algebra), satisfying a certain axiom, which can be explained in terms of the coproduct on the (Hopf) algebra of functions. One takes this axiom for the definition of an integral on an arbitrary Hopf algebra over an arbitrary field.

Since the pioneering work of Sweedler \cite{sweedler1}, integrals on Hopf algebras were studied by several authors \cite{sul,lin1,doi1,stefan1}. Among others, Sweedler proved the existence and uniqueness upto a constant of a (non-zero) integral on any finite-dimensional Hopf algebra. A theorem of Sullivan states that if an integral exists on a Hopf algebra then it is uniquely determinend up to a constant. Therefore we shall refer to the integral on a Hopf algebra when ever it exists. We shall also assume that integral means non-zero integral.

In representation theory of a compact group, one uses the Haar integral to deduce the semisimplicity of its representations. There is an analogue in comodule theory of Hopf algebras. If a Hopf algebra possesses an integral which does not vanish at the unit element, then it is cosemisimple, i.e., all its comodules are semisimple. However, the integral may vanish at the unit element, which is equivalent to the non-cosemisimplicity of the Hopf algebra. While there are many examples of non-cosemisimple finite-dimensional Hopf algebras, not so many infinite-dimensional non-cosemisimple Hopf algebras with integral are known. Moreover, a theorem of Sullivan \cite{sul} states that a commutative Hopf algebra over a field of characteristic zero  possesses an integral if an only if it is cosemisimple.

New examples of infinite-dimensional non-cosemisimple Hopf algebras with integral come from Lie supergroups and quantum (super)groups theory. In studying Haar measure on compact supergroups, F. Berezin found a remarkable fact that the Haar measure exists but the whole volume of the supergroup with respect to this measure may be zero (see e.g. \cite{berezin1}). In other words, the function algebra on a compact supergroup is a infinite-dimensional Hopf algebra with integral, which may not be cosemisimple.

In \cite{ph98b} the author showed that the Hopf algebras associated with certain (non-even) Hecke symmetries (i.e., quantum groups of type $A$) are non-cosemisimple infinite-dimensional Hopf algebras with integral.

In studying representations of simple Lie super-algebras of classical type, V. Kac found out that their irreducible representations split into two classes of typical and atypical representations (see, e.g., \cite{kac1}). It turns out that there is an analogous notion for simple comodules over a Hopf algebra with integrals and the integral provides a necessary and sufficient condition for a simple comodule to be ``typical'' (called ``splitting'' in this work). This is the main result of the first part of this work. In the second part we apply this result to study representations of quantum groups of type $A_{0|0}$, i.e., Hopf algebras associated to Hecke symmetries of birank $(1,1)$. Using the classification result we are able to classify the symmetries themselves.

The work is briefly divided into two parts. In order to reach to the main result of the first part (Theorem \ref{thm32}), we first recall some definitions and known facts on integrals on Hopf algebras (Section \ref{sect1}). Then we define a convolution product on a Hopf algebra by means of the integral making the Hopf algebra into a non-unital associated algebra and derive some auxiliary result for Section \ref{sect3} (Section \ref{sect2}). In Section \ref{sect3} we introduce the notion of splitting comodule, which means injective, projective simple comodule. In the terminology of V. Kac, splitting comodule splits in any comodule. We provide in Theorem \ref{thm32} a necessary and sufficient condition of a simple comodule to be splitting.

In the second part of the work, Section \ref{sect4}, we apply the result of the first part to Hopf algebras associated to Hecke symmetries of birank $(1,1)$, i.e., Hecke symmetries, the quantum exterior algebras associated to which have the Poincar\'e series equal to $(1+t)(1-t)^{-1}$. We show that simple comodules of these Hopf algebras can be labelled by pairs of integers $(k,l)$, where $(1,0)$ is the fundamental comodule, $(-1,0)$ is its dual, $(0,0)$ is the trivial comodule and the comodule labelled by $(k,l)$ is splitting iff $k+l\neq 0$. We show that the dimension of a simple comodule is 2 or 1 depending on whether it is splitting or not. Using this we able to classify the Hecke symmetries of birank $(1,1)$. In turns out that there are no other then those found by Manin \cite{manin2} and Takeuchi-Tambara \cite{tt}

\section{Co-Frobenius Coalgebras and Hopf algebras}\label{sect1}
We work over a field $\bK$. Every tensor product if not explicitly indicated means tensor product over $\bK$.

 Let $C$ be a coalgebra and $M$ be a right $C-$comodule, the coaction of $C$ on $M$ is denoted by $\rho$, $\rho:M\lora M\ot C$, $\rho(v)=v_0\ot v_1$. Let $C^*:=\Hom_\bK(C,\bK)$ be the dual of $C$. Then $C^*$ is an algebra, acting on $M$  from the left in the following way 
\bbs \phi\leftact v:=v_0\phi(v_1).\ees
Analogously, if $\lam:N\lora C\ot N$, $\lambda(v)=v_1\ot v_2$,  is a left $C-$comodule, then it is a right $C^*$-module through the action $v\rightact\phi:=\phi(v_1)\ot v_2$. 
%There, in fact, an equivalence between the categories of finite dimensional right $C-$comodules and of finite dimensional left $C-$comodules which appends to each right comodules structure on $M$ the left comodule structure on $M^*$: for $v\in M,\phi\in M^*$, $\lam(\phi)(v)=\phi(v_0)v_1$.

Thus, we have a functor from the category of right (left) $C-$comodules into the category of left (right) $C^*$-modules, which is full, faithfull and exact. The following statement is due to Doi.
\begin{proc}\label{doi} \cite{doi1}
{\it Let $M$ be a right (left) $C-$comodule, which is finite dimensional. Then $M$ is injective (projective) if and only if it is injective (projective) as left (right) ${C^*}$-module.}
\end{proc}

A $C^*$-module may not be a $C-$comodule by the above correspondence. Those $C^*$-module induced from  $C-$comodules are called rational modules. Each left $C^*$-module $M$ contains a unique maximal rational submodule, denoted by $_{\rm rat}M$. Analogously, for a right $C^*$-module $M$, its rational submodule will be denoted by $M_{\rm rat}.$

Let $M$ be a $C-$comodule. The map $\rho:M\lora M\ot C$ induces a map $M^*\ot M\lora C$, which can be considered as a coalgebra homomorphism or a morphism of $C-$comodules, where $C$ coacts on $M^*\ot M$ on the second tensor component. In the latter case, we shall use the notation $(M^*)\ot M$ to indicate that $C$ coacts only on $M$. The image of $M^*\ot M$ is called the coefficient space of $M$, denoted by $\Cf(M)$. 

Let $S$ be a simple (left or right) $C-$comodule. The fundamental theorem of comodule (saying that a finite generated comodule is finite dimensional) implies that $S$ is finite dimensional. Let $\D:=\End^C(S)$. Then, by Schur lemma, $\D$ is a division algebra over $\bK$ and $S$ is a vector space over $\D$. We have $\Cf(S)\cong V^*\ot_\D V$, as coalgebras \cite{green2}. Let $\{M_\alpha|\alpha\in\A\}$ be the set of all simple $C-$comodules. We define $\D_\alpha:=\End^C(M_\alph)$ and $m_\alpha:=\dim_{\D_\alpha}(M_\alpha)$, $d_\alpha:=\sqrt{\dim_\bK\D_\alpha}.$ Note that $d_\alpha$ are positive integers and
if $\bK$ is algebraically closed then $\D_\alpha=\bK$, $\forall \alpha$, i.e., $d_\alpha=1$.

By definition, the socle of a comodule $M$ is the sum of all its simple subcomodule. The sum is direct and is denoted by $\si(M)$. The injective hull (or cover) of $M$,  is by definition an injective comodule $\J(M)$ together with a morphism $M\lora \J(M)$ inducing an isomorphism $\si(M)\lora \si(\J(M))$. It is easy to see that the injective hull of a simple comodule, if it exists, is indecomposable. The following results are due to J. Green.

\begin{proc}\label{green}  \cite{green2}
(i) The injective hull of any comodule exists uniquely.
 
(ii) $C$ itself decomposes in to indecomposable injective subcomodules as follows
\bba\label{eqgreen} C\cong \bigoplus_{\alpha\in \A}\J(M_\alpha)^{\oplus m_\alpha}.\eea
 
(iii) If $C=\bigoplus_{\lam\in \L}N_\lambda$ is another decomposition then for each $\alpha\in \A$, the set $\{\lambda\in\L|N_\lam\cong\J(M_\alpha)\}$ contains exactly $m_\alpha$ elements.
\end{proc}

A bilinear form $b$ on  $C$ is called balanced if, for all $\phi\in C^*$,
\bbs b(x\rightact\phi,y)=b(x,\phi\leftact y).\ees
Balanced bilinear forms on $C$ are in 1-1 corespondence with right $C^*-$comodules homomorphism $r:C\lora C^*$ by the formula $r(x)(y)=b(x,y)$ and in 1-1 corespondece with left $C^*-$comodules homomorphism $r:C\lora C^*$ by the formula $l(x)(y)=b(y,x)$.

A coalgebra is called left (right) co-Frobenius if there exist  a left (right) monomorphism of $C^*$-modules $C\lora C^*$. The following results are due to B. Lin.
\begin{proc}\label{lin}\cite{lin1} If $C$ is a {\rm left} co-Frobenius coalgebra then:

(i) The  injective cover of every finite dimensional {\rm right} $C-$comodule is finite dimensional.

(ii) Every injective {\rm right} $C-$comodule is projective.

(iii) $_{\rm rat}C^*$  is dense in $C^*$.\end{proc}

In the next section we shall need the following result.
\begin{lem}\label{lem13} Let $M$ be a $C-$comodule of finite dimension. Then $M$ is projective (resp. injective) if and only if it is projective (resp. injective) in the category of finite dimensional $C-$comodules.\end{lem}
\proof It is sufficient to show the ``if'' part. Assume that $M$ is projective in the category of finite dimensional $C-$comodule. Consider a diagram
\bbs\bbar{rrl}&&M\\
&\exists\swarrow &\downarrow \pi\\
&N\stackrel{\nu}{\lora}&P\lora 0.\eear\ees
By replacing $P$ with $\Im(P)$ and $N$ with $\nu^{-1}(\Im(P))$, we can assume that $\pi$ is surjective, thus $P$ is finite dimensional. Let $\P$ be a basis of $P$ and $\N$ be a set of elements of $N$ such that $\nu(\N)=\P$. The submodule $N_1$ of $N$, generated by $\N$ is finite dimensional and we have $\nu(N_1)=P$. Hence, by assumption, there exists a morphism $\mu:M\lora N_1$: $\nu\circ \mu=\pi$.

Assume that $M$ is injective in the category of finite dimensional $C-$comodule. Consider the diagram
\bbs\bbar{rl}0\lora P&\stackrel{\nu}{\lora} N\\
\pi\downarrow &\swarrow\exists \\
M.&
\eear\ees

By replacing $P,N$ with $P/\Ker(\pi),N/\Ker(\pi)$, we can assume that $\pi$ injective. Hence $P$ has finite dimension. If $N/P$ is finite dimensional, we are done by the assumption on $M$. Otherwise, consider the set 
\bbs\A:=\left\{(N_\alpha,\mu_\alpha)|N_\alpha\supset P,\mu_\alpha:N_\lam\lora M, \mu_\lam\circ \nu=\pi\right\}.\ees
 Define an order on this set, setting $\alpha\prec \beta$ iff $N_\alpha\subset N_\beta$ and $\mu_\beta|_{N_\beta}=\mu_\alpha$. The chain condition is satisfied hence there exists a maximal element, say, $N_1$. Since any submodule of $N$, containing $P$ and having finited dimension is contained in $\A$, $N_1$ is strictly bigger then $P$. Were $N_1\neq N$, repeating the above process, we would get a submodule $N_1$, $N_2\succ N_1$, which contradicts the maximality of $N_1$. Hence $N_1=N$.\eee

Let now $H$ be a Hopf algebra. Then $\bK$ is a left (right) $H-$comodule be means of the unit map. A left (right) integral on $H$ is an $H-$comodule morphism $H\lora \bK$, where $H$ is considered as left (right) comodule on itself by means of the copoduct. Let $\int_r$ (resp. $\int_l$) denoted a left (resp. right) integral on $H$, then we have
\bba\label{lint} a_1\intl{a_2}=\intl{a}, \\
\label{rint} \intr{a_1}a_2=\intr{a},\eea
$\forall a\in H.$

We need the following information on the integrals.
\begin{proc}\label{sls}  Let $H$ be a Hopf algebra. The following conditions are equivalent

(i) $H$ possesses a left integral.
 
(ii) $H$ is left  co-Frobenius as a coalgebra.

(iii)  $H$ possesses a right integral.

(iv) $H$ is right co-Frobenius as a coalgebra.

(v)  The injective hull of every left comodule is finite dimension.

(vi)  The injective hull of every right comodule is finite dimension.

(vii) $H$ possesses a finite dimension injective left comodule.

(viii) $H$ possesses a finite dimension injective right comodule.

(ix) $H^*_{\rm rat}$  is dense in $H^*$.

(x) $_{\rm rat}H^*$  is dense in $H^*$.

(xi) Every injective left $H-$comodule is projective.

(xii) Every injective right $H-$comodule is projective.
\end{proc}
The first 6 conditions are due to Larson-Sweedler-Sullivan, the conditions (vii)-(xii) follow from \ref{lin}

Define bilinear form $b$: $b(x,y):=\intl{xS(y)}.$ Using the identity
\bba\label{first-id}
h_1\intl{h_2S(g)}=\intl{hS(g_1)}g_2\eea
which follows immediately form the definition of $\int_l$, we can easily show that $b$ is balanced. The folowing results are due to \c D. \c Stefan.
\begin{proc}\cite{stefan1} \label{stefan}Let $H$ be a Hopf algebra with integral. Then the following facts hold.

(i) The bilinear form $b$ is non-degenerate.

(ii) For any finite dimensional $H-$comodule,
$\dim_\bK(\Hom^H(H,M))=\dim_\bK M.$\end{proc}
In particular, we have

(i) the antipode is injective, and

(ii) there exists $h$ such that $\intl{S(h)}\neq 0$. Since $\int_l\circ S$ satisfies \rref{rint}, it is a right integral on $H$.

Assume for a moment that the field $\bK$ is algebraically closed. Let $R$ be the radical of $H$, i.e. $R=\oplus_\alpha \Cf(M_\alpha)$, where $\{M_\alpha,\alpha\in \A\}$ is the set of all simple left (or right) $H-$comodules. As we have seen in the previous subsection, $\Cf(M_\alpha)\cong M(m_\alpha)^*$, $m_\alpha:=\dim_{\bK}(M_\alpha)$. Fix idempotents $\{ e_{\alpha,i}|\alpha\in\A,1\leq i\leq m_\alpha\}$ of the algebras $\Mat(m_\alpha)$ -- the matrix ring of degree $m_\alpha$. They can be considered as linear functional on $R$ by defining
\bbs e_{\alpha,i}(\Cf(M_\beta))=0, \mbox{ whenever } \alpha\neq \beta.\ees
A theorem of Sweeder-Sullivan \cite{sul}, stating that there exists a coalgebra projection $H\lora R$, implies that $e_{\alpha,i}$ can be extended on the whole $H$ and that
\bbs H\rightact e_{\alpha,i}\quad (\mbox{ resp. } e_{\alpha,i}\leftact H) \mbox{ is a right (resp. left) $H-$comodules}.\ees
Consequently, we have a decompostion
\bbs H\cong \bigoplus_{\alpha\in\A,\atop{1\leq i\leq m_\alpha}}H\rightact e_{\alpha,i}\left(\mbox{ resp. }\cong\bigoplus_{\alpha\in\A,\atop{1\leq i\leq m_\alpha}}e_{\alpha,i}\leftact H\right)\ees
as right (resp. left) $H-$comodules. 

On the other hand, it is easy to see that $M_\alpha\subset H\rightact e_{\alpha,i}$ as right $H-$comodules. Thus, comparing with the decomposition in \ref{green}, we have:
\begin{proc}\label{pro}Assume that the field $\bK$ is algebraically closed. Then
\bba \J(M_\alpha)\cong H\rightact e_{\alpha,i}\label{sul-gr}\eea
as right $H-$comodules.\end{proc}

\section{The Convolution Product on $H$}\label{sect2}
We define a new product on $H$:
\bba\nonumber g*h:=h_1\intl{h_2S(g)}
=\intl{hS(g_1)}g_2\label{*product}\mbox{ (by \rref{first-id}). }\eea
Using \rref{first-id} we can easily check that $*$ is associative. $*$ is called the convolution product on $H$.
Denote  $\check{H}:=(H,*)$. Then $\check{H}$ is a (non-unital) algebra.
Let $V$ be a right $H-$comodule. Then $V$ is a left $H$-module by means of the action
\bbs h*v:=v_0\intl{v_1S(h)}.\ees
The verification again uses \rref{first-id}. Denote $\check V:=(V,*)$.

Let $f:V\lora W$ be a homomorphism of right $H-$comodules, i.e., $f(v)_0\ot f(v)_1=f(v_0)\ot v_1$. We have
\bbs h*f(v)=f(v)_0\intl{f(v)_1S(h)}
= f(v_0)\intl{v_1S(h)}
=f(h*v).\ees
Thus $f$ is a homomorphism of left $\check H$-modules. Conversely, if $f$ is a homomorphism $\check V\lora\check W$, then we have, for all $h\in H$,
\bbs f(v)_0\intl{f(v)_1S(h)}=f(v_0)\intl{v_1S(h)}.\ees
By the non-degeneracy of the integral (\ref{stefan}, (i)), we have
\bbs f(v)_0\ot{f(v)_1}=f(v_0)\ot{v_1},\ees
which means that $f$ is a homomorphism of right $H-$comodules. Thus we have
\bba\label{eqh-checkh}\Hom_{\check H}(\check V,\check W)=\Hom^H(V,W).\eea
In particular we have proved
\begin{lem}\label{checkh} If $M$ is a simple right $H-$comodule then $\check M$ is a simple left $\check H$-module.\end{lem}

Let now $V$ be a cyclic $\check H$-module, that is, there exists $\bar{v}\in V$, such  that $V$ is generated by $\bar v$. We want to define a coaction of $H$ on $V$. Let $v\in V$. Then there exists (not uniquely) $h\in \check H$, such that $v=h*\bar v$. Set $\hat\delta(v):=h_1*\bar v\ot h_2.$ We show that $\hat\delta$ is independent of the choice of $\bar v$ and $h$, and that it is in fact a coaction of $H$ on $V$.

The fact that $\hat\delta$ does not change when $\bar v$ is replaced by $\tilde v$ is represented by the equation
\bbs h_1*\bar v\ot h_2=(h*g)_1*\tilde v\ot (f*g)_2,\ees
where $g*\tilde v=\bar v$. This equation follows immediately form the definition of $*$. The independence on the choice of $h$ means that, whenever $h*\bar v=0$, we have $h_1*\bar v\ot h_2=0$. Indeed, we have
\bbs \intl{h_2S(g)}h_1*\bar v=(g*h)*\bar v=g*(h*\bar v)=0,\ees
for all $g\in H$. By the non-degeneracy of $\int$, we conclude that $h_1*\bar v\ot h_2=0$. 

The coassociativity and counitary of $\hat\delta$ are also checked directly using \rref{first-id}. Moreover, denoting by $\hat V$ the resulting $H-$comodule, we also have $\check{\hat V}\cong V$.

Since simple modules are cyclic, we have
\begin{lem}\label{cycle} Let $M$ be a simple $\check H$-module, then there exists a right $H-$comodule $\hat M$, such that $M\cong \check{\hat M}$.\end{lem}

We also need another action of $\check H$ on a right comodule $V$ of $H$, given by
\bbs h\circ v:=v_0\intl{hS(v_1)}.\ees
Indeed, we have
\bbas g\circ (h\circ v)&=& g\circ v_0\intl{hS(v_1)}= v_0\intl{gS(v_1)}\intl{hS(v_2)}\\
&=& v_0\intl{hS\left(\intl{gS(v_1)}v_2\right)}=v_0\intl{hS(g_1)}\intl{g_2S(v_1)}\\
&=&(g*h)\circ v.\eeas
It is again easy to checke that $\breve V:=(V,\circ)$ is a left $\check H$-module, and that, if $V$ is simple then $\breve V$ is simple. Lemma \ref{cycle} holds only in case the antipode is bijective. In fact, if $V$ is a cyclic $\check H$-module. Then we can define the following coaction of $H$ on $V$:
\bbs \dot\delta(v):=h_1\circ \bar v\ot S^{-2}(h_2),\ees
where $\bar v$ is a generating element and $h*\bar v=v$. Denote by $\dot V$ the comodule of $H$ induced from $V$, we have $\dot{\breve V}\cong V$.

Composing the operation $\breve{\  }$ and $\hat{\  }$ on a simple comodule $V$ we obtain a new simple comodule $\hat{\breve V}$, denoted by $V^\bullet$. The coaction of $H$ on $V^\bullet$ is given by
\bbs \delta^\bullet(v)=v_0\intl{h_1S(v_1)}\ot h_2,
\ees
with $h$ given by condition $v_0\intl{hS(v_1)}=v.$

Now, assume that $\bK$ is algebraically closed. Let $M_\alpha$ be a simple right $H-$comodule. Then $\breve M_\alpha$ is a simple left  $\check H-$comodule. The action of $\check H$ on $\breve M_\alpha$ induces a $\check H$-module homomorphism $\pi:\check H\lora\breve M_\alpha\ot (\breve M_\alpha^*)$, where $\breve M_\alpha\ot (\breve M_\alpha^*)\cong \breve M_\alpha^{\oplus \dim_\bK M_\alpha}$ as $\check H$-modules. The isomorphism \rref{eqh-checkh} shows that $\pi$ is a homomorphism of $H-$comodules. 

%Thus we showed that \bba\label{eqh-salpha} \Hom_{H}(H,M_\alpha)\neq 0.\eea

On the other hand, $H$ decomposes into the direct sum of its indecomposable injective subcomodules as in Lemma \ref{green}. For $h\in H\rightact e_{\beta,j}$, i.e., $h=e_{\beta,j}(g_1)g_2$ for some $g\in H$, and for $v\in\breve M_\alpha$, we have
\bbs h*v=v_0e_{\beta,j}(g_1)\intl{g_2S(v_1)}=v_0\intl{g(S(v_1)}e_{\beta,i}(v_2).\ees
Thus, if $\alpha\neq \beta$, $h*M_\alpha=0$, therefore, $\pi(h)=0$. Thus
\bba\label{eqjsahpha-s}\Hom^H(\J(M_\alpha),M^\bullet_\alpha)\neq 0.\eea
According to \ref{stefan},
\bba\label{eqstefan} \dim_\bK(\Hom^H(H,M^\bullet_\alpha))=\dim_\bK M^\bullet_\alpha.\eea
Since $H$ contains precisely $m_\alpha=\dim_\bK M_\alpha$ copies of $\J(M_\alpha)$, we conclude that
\bba &&\label{eqsalpha-sbeta1}
\dim_\bK\Hom^H(\J(M_\alpha),M^\bullet_\alpha)=1,\\
&&\label{eqsalpha-sbeta2}
\dim_\bK\Hom^H(\J(M_\beta),M^\bullet_\alpha)=0 \quad\mbox{if $\alpha\neq \beta$}.\eea

\begin{thm}\label{thm28} Let $H$ be a Hopf algebra with integral. Then for any simple comodules $M_\alpha,M_\beta$
\bbs \dim_\bK\Hom(\J(M_\alpha),M^\bullet_\beta)= \delta_\beta^\alpha d_\beta^2,\ees
where $d_\beta^2$  is the dimension over $\bK$ of $\D_\beta=\End^H(M_\beta)$. 
\end{thm}
\proof The case $\bK$ is algebraically closed is already proved. Assume that $\bK\neq \overline{\bK}$. Then $\D_\beta=\End^H(M_\beta)$ splits over $\overline\bK$: ${\D_\beta}\ot_\bK\overline\bK\cong \Mat_{\overline\bK}(d_\beta)$. For the extension $\overline{H}:=H\ot_\bK\overline{\bK}$, the comodule $\overline{M_\beta}:=M_\beta\ot_\bK\overline{\bK}$ decomposes into $d_\beta$ exemplars of the simple $\overline{H}-$comodule $M_\beta'$. Since $\overline{\J(M_\beta)}:=\J(M_\beta)\ot_\bK\overline{\bK}$ remains a direct summand of $\overline{H}$, it is an injective $\overline{H}-$comodule. Therefore $\overline{\J(M_\beta)}$ is a direct sum of $d_\beta$ exemplars of $\J(M_\beta')$. Since, for $\alpha\neq \beta$,
\bbs \Hom_{\overline H}(\J(M_\alpha'),{M^\bullet_\beta}')=0,\ees
we have 
\bbs 
\Hom_{H}(\J(M_\alpha),M^\bullet_\beta)=0.\ees
Therefore, by virtue of Equation \rref{eqstefan} (which is valid on any field),
\bbs 
\dim_\bK\Hom_{H}(\J(M_\alpha)^{\oplus m_\alpha},M_\alpha)=\dim_\bK M_\alpha= d_\alpha^2 m_\alpha.\ees
Consequently
\bbs\dim_\bK \Hom_{H}(\J(M_\alpha),M_\alpha)= d_\alpha^2.\ees
\eee

 Let $M$ be a finite dimensional right $H-$comodule then $M^*:=\Hom_\bK(M,\bK)$ is also a right comodule with the coaction give by the equation
\bbs \rho(\phi)(x):=\phi_0(x)\phi_1=\phi(x_0)S(x_1),\quad x\in M,\phi\in M^*.\ees
The map $\ev:M^*\ot M\lora \bK$, $\phi\ot x\loma \phi(x)$ is a morphism of $H-$comodules. The pair $(M^*,\ev)$ is called left dual to $M$, it is defined uniquely upto isomorphism.
There exists a monomorphism $\db:\bK\lora S\ot S^*$, defined by the conditions $(\ev\ot\id_{S^*})(\id_{S^*}\ot \db)=\id_{S^*}$ and $(\id_s\ot \ev)(\db\ot\id_S)=\id_S$, which is also a comodule morphism. Dually, $(M,\ev)$ is called the right dual to $M^*$.

Thus, we see that the left dual to a finite dimensional comodule always exists. If the antipoded is bijective then the right dual to any finite dimensional comodule also exists. We shall need the following isomorphism, given by manipulating the morphism $\ev$ and $\db$: for any finite dimensional comodule $N$,
\bba \label{eq7}& \Hom^H(M\ot N,P)\cong\Hom^H(M,P\ot N^*)&\\
&\Hom^H(M,N\ot P)\cong\Hom^H(N^*\ot M,P).&\label{eq8}\eea
As an immediate corollary of Lemma \ref{lem13} and these equations, we have
\begin{lem}\label{lem29}Let $M$ be a finite dimensional comodule. Then we have:

(i)  If $M^*$ is projective (resp. injective) then $M$ is injective (resp. projective).

(ii) If the antipode is bijective, then $M$ is injective (resp. projective) iff $M^*$ is projective (resp. injective).\end{lem}
\proof Equations \rref{eq7} and \rref{eq8} imply
\bba\label{eq9} \Hom^H(M^*,N^*)\cong\Hom^H(N,M).\eea
Thus, if $M$ is projective (resp. injective) then $M^*$ is injective (resp. projective) in the category of finite dimensional comodules.\eee
\begin{cor}\label{cor-jhull} Let $H$ be a Hopf algebra wiht integral. Assume that the antipode is injective. Then for any simple comodule $M$
\bbs \J((M^\bullet)^*)\cong \J(M)^*.\ees
\end{cor}

\begin{pro}\label{vbullet} Assume that the Hopf algebra $H$ as a left-right integral. Then $V^\bullet\cong V^{**}$. If $H$ is moreover coquasitriangular then $V^\bullet\cong V$ and, consequently, $\J(M^*)\cong\J(M)^*$.\end{pro}
\proof 
Assume that $\int$ is a left-right integral. Thus we can define $V^\bullet$ as above. We want to show that
\bbs v_0\intl{h_1S(v_1)}\ot h_2=c\ot v_0\ot S^2(v_1),\ees
for certain constant $c$, depending only on $\int$. By the non-degeneracy of integral, this equation is equivalent to
\bbs v_0\int(v_1S(h))\int( h_2S(g))=c\ot v_0\int(S^2(v_1)S(g)).\ees
We have
\bbas&& \lefteqn{v_0\int(h_1S(v_1))\int( h_2S(g))}\\
&&=v_0\int((g*h)S(v_1))=(g*h)\circ v=g*(h\circ v)=g\circ v=v_0\int(gS(v_1)).
\eeas
By the uniqueness of integral, we can choose $c$ such that 
$c\ot\int(S^2(v_1)S(g))=\int(gS(v_1))$.

If $H$ is coquasitriangular then $V\cong V^{**}$. \eee
\begin{rem}\rm If the left and the right integrals do not coincide then in general, $M^\bullet\not\cong M^{**}$. An example is Sweedler's Hopf algebras, see, e.g., \cite{schmuedgen1}.\end{rem}

\section{Splitting Comodules}\label{sect3}
Let $S$ be a simple comodule over $H$. $S$ is called splitting comodule, or typical comodule, if $S=\J(S)$. Since $\J(S)$ is injective and hence projective, we see that $S$ splits in any comodule. This explain the name splitting. The name typical was used by V. Kac for modules over a Lie superalgebra \cite{kac1}.

 By virtue of conditions in \ref{sls}, if a Hopf algebra possesses a splitting comodule then it possesses a non-zero integral. The converse statement is not true. The aim of this section is to give a criteria for a simple comodule to be typical. 

Let $M$ is a right $H-$comodule then the coaction of $H$ on $M^{**}$ -- the double left dual to $M$ is given by (identifying $M^{**}$ with $M$ as vector spaces)
\bbs \rho_{M^{**}}(v)=x_0\ot S^2(v_1).\ees

\begin{lem}\label{lem22} Let $M$ be a simple $H-$comodule. Then $M$ is splitting iff $M^*$ is splitting.\end{lem}
\proof
 Assume that $M$ is a typical. Then $M$ is injective. By a theorem of Doi \cite{doi1}, $M^*\ot M$ is injective, too. By definition of $M^*$, we have an epimorphism $\ev:M^*\ot M\lora \bK$. Among indecomposable injective subcomodule of $M^*\ot M$ there exists one, say $J$, such that the restriction of $\ev$ on $J$ is not zero. On the other hand, since $J$ is indecomposable and injective, it should appear in the decomposition \rref{eqgreen}, and by Theorem \ref{thm28}, the only comodule with this property is $\J(\bK)$, the injective hull of $\bK$. Thus we show that $M^*\ot M$ contains $\J(\bK)$ as a subcomodule, consequently, it contains $\bK$ as subcomodule, i.e. $\Hom^H(\bK,M^{*}\ot M)\neq 0$. According to \rref{eq8}, we have
\bbs \Hom^H(M^{**},M)\neq 0.\ees
Therefore $M^{**}\cong M$ and hence is splitting. Consequently $M^*$ is also splitting, by Lemma \ref{lem29}.

Assume now that $M^*$ is splitting. The discussion above shows that $M^{**}$ is also splitting and $M^{**}\cong M^{****}$. Since
\bbs \Hom^H(M^*,N^*)\cong\Hom^H(N,M),\ees
we conclude that $M\cong M^{**}$. Thus $M$ is splitting.\eee

\begin{thm}\label{thm32}
Let $M$ be a simple right $H-$comodule.
  Then $M$ is splitting if and only if the bilinear form $c$, $c(x,y)=\intr{yS(x)}$, is not identically zero on $\Cf(M).$ In this case, $c$ is also non-degenerate on $\Cf(M)$.\end{thm}
\proof 
``if''.

For each $g\in H$, define a linear functional $\phi_g\in\H^*$: $\phi_g(h):=\intr{hS(g)}$. By assumption, there exist $g\in\Cf(M)$ such that $\phi_g$ is not identically zero on $\Cf(M)$. Since $M$ is simple, the right coideal generated by $g$ is isomorphic to $M$. Define a linear map $\eta=\eta_g$:
\bba\label{eqeta}
\eta:\Cf(M)\lora M^{**},\quad h\loma g_1\intr{hS(g_2)}, h\in M.\eea
Since $\va(\eta(h))=\phi_q(h)$, $\eta$ is not trivial.

We have the following identity, which is an immediate consequence of \rref{rint} and the injectivity of the antipode
\bba\label{second-id}
\intr{h_1S(g)}h_2=S^2(g_1)\intr{hS(g_2)}.\eea
It follows form \rref{second-id} that $\eta$ is a morphism of $H-$comodules. Since $\eta$ is non-trivial on $\Cf(M)$, which is a direct sum of copies of $M$, we conclude that $\eta$ should induce a morphism $M\lora M^{**}$, which is non-trival. Since $M$ is simple and $\dim_\bK M=\dim_\bK M^{**}$, this morphism is an isomorphism. As a consequence, $\Cf(M^{**})=S(\Cf(M))=\Cf(M)$ and $\phi_g$ is $*$-invertible on $\Cf(M)$. Let $\psi$ be the $*$-inverse to $\phi$, define on $\Cf(M)$, thus $\phi(h_1)\psi(h_2)=\psi(h_1)\phi(h_2)=\va(h).$

Let now $M\hookrightarrow N$ be an inclusion of $H-$comodules. Let $f:N\lora M$ be a linear projection on $M$. We define a new map $F:N\lora M$ as follows,
\bbs F(v):=f(v_0)_0\psi(f(v_0)_1)\intr{v_1S(f(v_0)_2)}.\ees
$F$ is well defined by the assumption that $\Im(f)=M$, which implies $f(v_0)_1\in \Cf(M).$ For $v\in M$, $f(v)=v$, hence
\bbs F(v)=v_0\psi(v_1)\intr{v_3S(v_2)}=v_0\psi(v_1)q(v_2)=v.\ees
Thus, $F$ is again a projection of $M$. If we show that $F$ is a morphim of $H-$comodule, then we will be done.

By defintion of $\phi$, we can consider $F$ as a composition of the map $g:N\lora M^{**}$:
\bbs g(v)=f(v_0)_0\intr{v_1S(f(v_0)_1},\ees
and the morphism $\eta^{-1}:M^{**}\lora M$. Thus, it is sufficient to show that $g$ is a morphism of $H-$comodules, which means
\bbs f(v_0)_0\intr{v_1S(f(v_0)_1}\ot v_2=f(v_0)_0\intr{v_1S(f(v_0)_2}\ot S^2(f(v_0)_1).\ees
We have, according to \rref{second-id},
\bbas \mbox{ the left-hand side }&=& f(v_0)_0\ot S^2(f(v_0)_1)\intr{v_1S(f(v_0)_2)}\\
&=&\mbox{ the right-hand side }.\eeas
Therefore, $F$ is a morphism of $H-$comodules, consequently, $M$ is injective and $\J(M)=M$.

``only if''\ \  Assume now that $M$ is splitting, then, by Lemma \ref{lem22},  $M^*$ is also splitting, hence $M\ot M^*$ is injective. By definition of $M^*$, there exists a monomorphism $\db:\bK\lora M\ot M^*$, inducing a monomorphism $\J(\bK)\hookrightarrow M\ot M^*$. The latter inclusion induces the following inclusion
\bbs\\J(\bK)\subset\Cf(\J(\bK))\hookrightarrow \Cf(M\ot M^*)=\Cf(M)\ot\Cf(M^*).\ees
Since the right integral does not vanish identically on $\J(\bK)$ (by \ref{thm28}), we conlude that the set $\intr{a^i_jS(a^k_l)}$ is not identically zero, as $a^i_jS(a^k_l)$ span $\Cf(M\ot M^*)$. \eee

\section{Simple Representations of  Quantum Groups of Type $A_{0|0}$} \label{sect4}

Let $V$ be a finite dimensional vector space over $\bK$, a field of characteristic zero. An operator $R:V\ot V\lora V\ot V$ is called a Hecke symmetry if $R$ satisfies the Yang-Baxter equation 
$$ (R\ot \id_V)(\id_V\ot R)(R\ot \id_V)=(\id_V\ot R)(R\ot \id_V)(\id_V\ot R),$$
the Hecke  equation
$$(R-\id)(R-q\cdot \id)=0,\quad q\neq 0$$
and is closed, that is, the operator $P:V^*\ot V\lora V\ot V^*$,  half dual (half-adjoint) to $R$ -- given by 
$$ P=(\ev_V\ot\id_{V\ot V^*})(\id_{V^*}\ot R\ot\id_{V^*})(\id_{V^*\ot V}\ot\db_V)$$
 is invertible. We shall also assume that $q$ is not a root of unity of degree greater than 1.

Being given a Hecke symmetry, one can define the associated quantum symmetric and anti-symmetric tensor algebras as factor algebras of the tensor algebra over $V$ by the ideal, generated by $\Im(R-q\dot\id)$ and $\Im(R+\id)$, respectively. It is shown that the Poincar\'e series of these algebras, i.e., the formal sums with coefficients being dimensions of homongeneous components of these algebras, are rational functions \cite{ph97c} having negative roots and positive poles.

A Hecke symmetry $R$ is said to have birank $(1,1)$ if the Poincar\'e series of the associated quantum symmetric tensor algebra has one pole and one root, i.e., is of the form $(1+at)(1-bt)^{-1}, a,b>0$.

The quantum group (quantum semi-group) associated to $R$  is defined to be the Hopf algebra (bialgebra) universally coacting on the mentioned above quantum tensor algebras \cite{manin1}.  They are denoted by $H$ and $E$, respectively. If $R$ has birank $(1,1)$, the associated quantum group is called quantum group of type $A_{0|0}$. 
Simple $E-$comodules can be labelled by hook-partitions of the form $(m,1^n)$, $m\geq 1,n\geq 0$ and the trivial partion $(0)$. For simplicity we shall use the pair $(m,n)$ to denote the partion $(m,1^n)$ and the pair $(0,0)$ to denote the trivial partition. The endomorphims ring of a simple $E-$comodule is isomorphic to $\bK$. On the other hand, simple $E-$comodules are also simple as $H-$comodules with the natural action induced from the inclusion $E\lora H$. The reader is referred to \cite{ph97b} for detail.

From now on, for simplicity we shall use a dot $\cdot$ to denote the tensor product and a plus $+$ to denote the direct sum, thus $V^n$ will means $V^{\ot n}$ and $n\cdot V$ will means $V^{\oplus n}$. We  shall also use the equal sign $=$ to denote an isomorphism.

\def\ot{\cdot}
\def\oplus{+}
Simple  $E-$comodules  $\I{m}{n}$, associated to pairs $(m,n), m\geq 1, n\geq 0$, are given by the following rule. $\I{n}{0}=S_n$ is the $n$-th component of the quantum symmetric tensor over $V$, $\I{1}{n-1}=\Lambda_n$ is the $n$-th component of the quantum anti-symmetric tensor over $V$, $n\geq1$, $\I{0}{0}:=\bK$,
\bba\label{neq1}
\I{p}{q}\ot\I{m}{n}=\I{m+p}{n+q}\op \I{m+p-1}{n+q+1},\eea
for $m,p\geq 1, n,q\geq 0$. Particularly, we have 
$V=\I10$, $V^*=\I{-1}0$.

Our aim is to associate to each pair $(m,n)$ of integers a simple $H-$comodule and show that they furnish all simple $H-$comodules.

According to a result of \cite{ph98b}, if $\rank_qR=0$, where $\rank_qR$ is the full trace of the half-dual opertor $P$, then $H_R$ possesses an integral. Thus, in order to apply the results of the previous section, we have to show that $\rank_qR=0$. To do this we consider the Koszul complex of the second type introduced by Manin \cite{manin3} (see also \cite{gur1,ph98a,ls}). It is shown that, if $\rk_q R\neq -[k-l]_q$, then the complex
\bbs \cdots\lora\Lambda_k\ot S_l^*\stackrel{d_{k,l}}{\lora}\Lambda_{k+1}\ot S_{l+1}^*\stackrel{d_{k+1,l+1}}{\lora}\Lambda_{k+2}\ot S_{l+2}^*\lora\cdots\ees
with the differential induced from the dual basis map $\ev:\bK\lora V\ot V^*$, is exact \cite{gur1,ls}.

Notice that, according to \rref{neq1}, \rref{eq7} and \rref{eq8}, for $m,p\geq 1$ and $n,q\geq 0$,
\bbas&&\Hom^H(\I{p}{q},\I{m+p}{n+q}\ot \I{m}{n}^*)=\bK\\
&& \End^H(\I{m+p}{n+q}\ot \I{m}{n}^*)=\End^H(\I{m+p}{n+q}\ot \I{m}{n})=2\cdot \bK.
\eeas
Therefore, denoting $\I{-m}{-n}:=\I{m}{n}^*$, we have, for $p,q\geq 1$, $m\geq 1,n\geq 0$,
\bba\label{neq2}\I{p}{q}\op\I{p+1}{q-1}\subset \I{m+p}{n+q}\ot\I{-m}{-n},\eea
and $\I{m+p}{n+q}\ot\I{-m}{-n}$ does not contain any other simple $E-$comodule.

Assume that $\rk_qR\neq 0$. Then the complex
\bba\label{neq3}0{\sora} \bK\stackrel{d_{00}}{\sora}\I{1}{0}\ot\I{-1}{0}\stackrel{d_{11}}{\sora}\I{1}{1}\ot\I{-2}{0}\stackrel{d_{11}}{\sora}\cdots\eea
is exact. We have, for $n>m\geq 1$,
\bba\nonumber \I{n}{1}\ot\I{1}{m-1}\ot\I{-m}{0}&=&(\I{n+1}{m}\op\I{n}{m+1})\ot\I{-m}{0}\\
&\supset&2\cdot\I{n-m+1}{m}\op\I{n-m}{m+1}\op\I{n-m+2}{m-1}\label{neq30}.\eea
Thus, multiplying  \rref{neq3} with $\I{n}{1}$, we have a diagram
\bba\label{neq31}\bbar{rcccccc}
\I{n}{1}&&2\cdot\I{n}{1}\op\I{n-1}{2}\op\I{n+1}{0}&&
2\cdot \I{n-1}{2}\op\I{n-2}{3}\op\I{n}{1}&\dots\\
\|&&\bigcap&&\bigcap&&\\
0\sora\I{n}{1}&\sora&\I{n}{1}\ot\I{1}{0}\ot\I{-1}{0}&\sora&\I{n}{1}\ot\I{1}{1}\ot\I{-2}{0}&\sora\cdots
\eear\eea
The exactness of the lower complex and the remark following \rref{neq2} imply that $\I{n+1}{0}=0$, contradiction. Thus we have $\rk_qR=0$. As a consequence, the Hopf algebra $H$ possesses an integral and the formula for the integral in \cite{ph98b} implies that $\I{m}{n}$, and hence $I_{-m,-n},m\geq 1,n\geq0$,  are all splitting, except for $\I{0}{0}=\bK$. Therefore, by means of the two isomorphism preceeding \rref{neq2}, the inclusion in \rref{neq2} is in fact an isomorphism: for $p,q,m\geq 1,n\geq 0,$
\bba\label{neq2'}\I{m+p}{n+q}\ot\I{-m}{-n}=\I{p}{q}\op\I{p+1}{q-1} ..\eea

The next step is to define the comodules $\I{-1}{1}$ and $\I{1}{-1}$.

Consider the sequence \rref{neq3}.
Since $\rk_qR=0$, $\I{1}{0}\ot\I{-1}{0} =V\ot V^*$ contains two exemplars of $\bK$ in its composition series but only one as subcomodule. Let $M:=(V\ot V^*)/\bK$, then the map $V\ot V^*\lora \bK$ factorizes through $\bK$ to a map $M\lora\bK$. Dualizing this we get a sequence $\bK\lora M^*\lora V\ot V^*$. Since 
$$\Hom^H(\bK,\I{1}{1}\ot\I{-2}{0})\cong\Hom^H(\I{2}{0},\I{1}{1})=0,$$
that is $\bK$ cannot be a subcomodule of $\I{-2}{0}$, $\Im d_{1,1}\neq \bK$. From \rref{neq31}, we see that $\Ker d_{2,2}{{\supset}\atop{\neq}}\bK$. Let $N:=(V\ot V^*)/\Ker d_{2,2}$. Then $N$ is a factor comodule of $M$, which is different from $\bK$ and $M$. Therefore $V\ot V^*$ contains at least 4 simple comodules in its composition series. It cannot be larger then 4, as on the left hand side of \rref{neq30}, there are 4 simple comodules.
Denote by $A$ and $B$ the two simple subcomodules, that are different from $\bK$. Since $V\ot V^*$ is self dual, either $A$ and $B$ are self dual or $B$ is dual to $A$.

Using \rref{neq2'}, we have
\bbs \I{2}{1}\ot \I{-1}{0}\ot\I{1}{0}=2\ot\I{2}{1}+\I{1}{2}+\I{3}{0}.\ees
Thus, we can assume that $\I{2}{1}\ot A=\I{1}{2}$ and $\I{2}{1}\ot B=\I{3}{0}.$ Using induction we can easily show that
\bba\label{neq32}\I{m}{n}\ot A=\I{m-1}{n+1}&& \I{m}{n}\ot B=\I{m+1}{n-1},\eea
for all $m\geq 2,n\geq 1$. Using the fact, that $M\ot M^*$ contains $\bK$ as a subcomodule, for any comodule $M$, we deduce that $A= B^*$ and $A\ot B=\bK$. Thus
\bba\label{neq4} &&\I{1}{-1}\ot \I{-1}{1}=\bK\\
\label{neq5} && \I mn\ot \I 1{-1}=\I{m+1}{n-1}\\
\label{neq6} && \I mn\ot \I {-1}1=\I{m-1}{n+1},\eea
for all $m\geq 2,n\geq 1$. Dualizing these equalities, we obtain
\bba\label{neq7} &&\I{-m}{-n}\ot \I{1}{-1}=\I{-m+1}{-n-1}\\
\label{neq8} && \I{-m}{-n}\ot \I{-1}1=\I{-m-1}{-n+1},\eea
for all $m\geq 2,n\geq 1$. 

Consider now
\bbas \I11\ot \I10\ot \I{-1}0&=&(\I21+\I12)\ot \I{-1}0\\
&=&\I11+\I20+\I12\ot \I{-1}0.\eeas
The left-hand side contains simple comodule $\I11\ot\I1{-1}$ and $\I11\ot \I{-1}1$. We therefore conclude that
\bbs \I11\ot \I1{-1}=\I20,\ees
and thus
\bbs \I20\ot \I{-1}1=\I11.\ees

We are now at the stage to associate to each pair $(m,n)$ of integers a simple comodule $\I mn$.
Note that for $m\neq 1,n\geq 0$ or $m\leq\-1,n\leq 0$, we have already define $\I mn$. We call $s(m,n):=m+n$ the total degree of the pair $(m,n)$. Thus, there can be three possiblities: $s(m,n)>0; <0$ or $=0$. If $s(m,n)=0$, i.e., $m=-n$, set
\bbs \I m{-m}:=\I1{-1}^m.\ees
If $s(m,n)\neq 0$, set
\bbs \I mn:=\I{m+n}0\ot \I{-n}n.\ees

Using (\ref{neq4}-\ref{neq8}), it is easy to see that the definition is compatible with the predefined comodules and that these comodules are all simple. We want to find the formula for the tensor product of these comodules and deduce from this formula that these comodules furnish all simple $H-$comodule.

Let $(m,n)$ and $(p,q)$ be pairs of integers. Our aim is to decompose $\I mn\ot \I pq$. The main role here plays the total degree. There can be three possibilities
\begin{enumerate}\item either $m+n$ or $p+q$ is equal to zero;
\item $m+n$ and $p+q$ are both different from zero but their sum is zero;
\item $m+n$ and $p+q$  and $m+n+p+q$ are all different from zero.\end{enumerate}

1. If $m+n=0$ then $\I mn=\I1{-1}^m$, hence
\bba\label{neq9} \I m{-m}\ot\I pq=\I{p+m}{q-m}.\eea

2. If $m+n+p+q=0$ and $m+n\neq 0$, using (\ref{neq4}-\ref{neq8}), we can assume $n=p=0$. Thus $m=-p$. We claim that 
\bba\label{neq10} \I m0\ot \I{-m}0=\I 01+\I10.\eea
Indeed, $\I m0^*=\I{-m}0$, hence $\I m0\ot \I{-m}0$ contains $\bK$ as subcomodule. More over, this comodule is injective. On the other hand, $\I10\ot \I{-1}0$ is the injective envelope of $\bK$, therefore (cf. \cite{green2}), is a subcomodule of $\I m0\ot \I{-m}0$, $\forall m\geq 0$. Multiplying these comodules with $\I m1$ we get the same comodule. Whence we conclude \rref{neq10}.

3. If $m+n$, $p+q$, $m+n+p+q$ are all non-zero, dualizing if necessary, we can assume $m+n+p+q>0$. Using (\ref{neq4}-\ref{neq8}), we can assume $n=q=0$. Assume $m>p$, thus $m>0$. One is led to computing $\I m0\ot \I p0$. If $p>0$, the formula is already known (cf. \ref{neq1}-\ref{neq2}). Assume $p< 0$ and set $k=-p$, then $k>0$ and $m> k$. We consider two case: $m-k\geq 2$ and $m-k=1$. If $m-k\geq 2$, then, according to (\ref{neq4}-\ref{neq8}),
\bba \I m0\ot \I{-k}0&=& \I1{-1}\ot \I{m-1}1\ot \I{-k}0\nonumber\\
&=& \I1{-1}\ot(\I{m-k-1}1+\I{m-k}0)\nonumber\\
&=&\I{m-k}0+\I{m-k+1}{-1}.\label{neq11}\eea
In the case $m-k=1$, we show that
\bba\label{neq12} \I m0\ot \I{-m}0=\I10+\I2{-1}.\eea
We have
\bbas \Hom(\I10,\I m0\ot\I{-m+1}0)&=& \Hom(\I10\ot \I{m-1}0,\I m0)\\
&=& \bK.\eeas
Remember that $\I2{-1}=\I1{-1}\ot \I10$ and  that $\I2{-1}$ is also simple. Since
\bbas \Hom(\I m0\ot \I{-m+1}0,\I01)&=& \Hom(\I{-1}1\ot \I {m}0,\I{m-1}0\ot \I10)\\
&=& \Hom(\I{m-1},\I{m-1}1\ot \I10)\\
&=& \bK.\eeas
Thus $\I m0\ot \I{-m+1}0$ contains $\I10$ and $\I2{-1}$ as subcomodules. On the other hand, multiplying both sides of \rref{neq12} with $\I{m+2}1,$ we get an equality. Therefore \rref{neq12} is proven. We summarize the results obtained in a theorem.
\begin{thm}\label{classification} Simple representation of a quantum group of type $A_{0|0}$ are classified by pairs $(m,n)$ of integers with the following properties:
\begin{enumerate} \item $\I m0$ is the $n$-th symmetric tensor, $\I{1}{n-1}$ is the $n$-th anti-symmetric tensor, $\I00=\bK$, $\I mn^*=\I{-m}{-n}$, $\I1{-1}$ is the super determinant. $\I mn$ is splitting iff $m+n\neq 0$.
\item We have the following rule for tensor product of simple comodules.
\begin{enumerate}\item for any integers $m,n$,
$$ \I mn\cdot \I{-1}1=\I{m-1}{n+1},$$
\item for any $m>n>0$,
\bbas &\I m0\cdot \I n0=\I{m+n}0+\I{m+n-1}1&\\
& \I m0\cdot \I{-n}0=\I{m-n}0+\I{m-n+1}{-1}&.\eeas
\item for $m\neq 0$,
$ \I m0\cdot \I{-m}0=\I 10\cdot\I{-1}0$, this comodule is injective and indecomposable. It contains two exemplars of $\bK$ and the comodules $\I1{-1}$ and $\I{-1}1$ in its decomposition series.\end{enumerate}\end{enumerate}\end{thm}

The classification obtained above also allowes us to classify Hecke symmetries of birank $(1,1)$.
The crucial point here is to compute the dimension of simple comodules. Since $\I 1{-1}\cdot \I{-1}1=\I00=\bK$, $\I1{-1}$ is one-dimensional. On the other hand, assuming that the Poincar\'e series of the quantum anti-symmtric algebra $\wedge$ is $(1+at)(1-bt)^{-1}$ with $a,b>0$, we can compute the dimension of polynomial comodules $\I mn$, $m\geq 1,n\geq 0$,
$$\dim_\bK\I mn=a^mb^n+a^{m-1}b^{n+1}.$$
According to \rref{neq5}, we have $a=b$. On the other hand, computing the dimension of $\I10\cdot \I{-1}0$ in two ways we obtain $a+b=2$. Therefore $a=b=2$, that is $\dim_\bK V=2$. That means, a Hecke symmetry of birank $(1,1)$ should be defined on a vector space of dimension $2$. There are only two families of such operators. The first one is two-parementeric, found by Manin \cite{manin2}, the second one is one-paramentric, found by Tambara-Takeuchi \cite{tt}.

\begin{center}\bf Acknowledgment\end{center}

The work was done during the author's stay at the Max-Planck-Institut f\"ur Mathematik, Bonn, Germany.


\begin{thebibliography}{10}

\bibitem{berezin1}
Felix~Alexandrovich Berezin.
\newblock {\em Introduction to superanalysis}.
\newblock D. Reidel Publishing Co., Dordrecht, 1987.

\bibitem{doi1}
Yukio Doi.
\newblock Homological coalgebra.
\newblock {\em J. Math. Soc. Japan}, 33(1):31--50, 1981.

\bibitem{green2}
J.A. Green.
\newblock {Locally Finite Representations}.
\newblock {\em Journal of Algebra}, 41:137--171, 1976.

\bibitem{gur1}
D.I. Gurevich.
\newblock { Algebraic Aspects of the Quantum Yang-Baxter Equation}.
\newblock {\em Leningrad Math. Journal}, 2(4):801--828, 1991.

\bibitem{ph98a}
Phung~Ho Hai.
\newblock {Hecke Symmetries}.
\newblock {\em J. of Pure and Appl. Algebra}.
\newblock to appear.

\bibitem{ph97c}
Phung~Ho Hai.
\newblock {Poincar\'e Series of Quantum Spaces Associated to Hecke Operators}.
\newblock {\em Acta Math. Vietnam., to appear}.
\newblock Available at {\tt xxx.lanl.gov/dvi/q-alg/9711020}.

\bibitem{ph97b}
Phung~Ho Hai.
\newblock {On Matrix Quantum Groups of Type $A_n$}.
\newblock {\em Preprint ICTP/97101, {\tt q-alg/9708007}}, 1997.

\bibitem{ph98b}
Phung~Ho Hai.
\newblock {The integral on quantum super groups of type $A_{r|s}$ }.
\newblock {\em {\sl Preprint MPIM 1998/127}}, 1998.

\bibitem{kac1}
V.~G. Kac.
\newblock Characters of typical representations of classical {L}ie
  superalgebras.
\newblock {\em Comm. Algebra}, 5(8):889--897, 1977.

\bibitem{schmuedgen1}
Anatoli Klimyk and Konrad Schm{\"u}dgen.
\newblock {\em Quantum groups and their representations}.
\newblock Springer-Verlag, Berlin, 1997.

\bibitem{lin1}
Bertrand I-peng Lin.
\newblock Semiperfect coalgebras.
\newblock {\em J. Algebra}, 49(2):357--373, 1977.

\bibitem{ls}
V.V. {L}yubashenko and A.~Sudbery.
\newblock {Quantum Super Groups of {GL}$(n|m)$ Type: Differential Forms,
  {K}oszul Complexes and {B}erezinians}.
\newblock {\em Duke Math. Journal}, 90:1--62, 1997.

\bibitem{manin3}
Yu.I. {M}anin.
\newblock {\em Gauge Field Theory and Complex Geometry}.
\newblock Springer-{V}erlag, 1988.

\bibitem{manin1}
Yu.I. {M}anin.
\newblock {\em {Quantum Groups and {N}on-commutative {G}eometry}}.
\newblock GRM, Univ. de Montreal, 1988.

\bibitem{manin2}
Yu.I. {M}anin.
\newblock {Multiparametric Quantum Deformation of the {G}eneral Linear
  Supergroups}.
\newblock {\em Comm. Math. Phys.}, 123:163--175, 1989.

\bibitem{stefan1}
Drago{\c{s}} {\c{S}}tefan.
\newblock The uniqueness of integrals (a homological approach).
\newblock {\em Comm. Algebra}, 23(5):1657--1662, 1995.

\bibitem{sul}
J.B. Sullivan.
\newblock {The Uniqueness of Integral for {H}opf Algebras and Some Existence
  Theorems of Integrals for Commutative Hopf Algebras}.
\newblock {\em Journal of Algebra}, 19:426--440, 1971.

\bibitem{sweedler1}
Moss~Eisenberg Sweedler.
\newblock Integrals for {H}opf algebras.
\newblock {\em Ann. of Math. (2)}, 89:323--335, 1969.

\bibitem{tt}
M.~{T}akeuchi and D.~{T}ambara.
\newblock A new one-parameter family of $2\times 2$ quantum matrices.
\newblock {\em Hokkaido Math. Journal}, XXI(3):409--419, 1992.
\newblock See also Proc. Japan. Acad., 8(8), 1991.

\end{thebibliography}
\end{document}